# Cowles commission structural equation approach in light of nonstationary time series analysis

## Cheng Hsiao[1]

*University of Southern California*

**Abstract:** We review the advancement of nonstationary time series analysis from the perspective of Cowles Commission structural equation approach. We argue that despite the rich repertoire nonstationary time series analysis provides to analyze how do variables respond dynamically to shocks through the decomposition of a dynamic system into long-run and short-run relations, nonstationarity does not invalid the classical concerns of structural equation modeling — identification and simultaneity bias. The same rank condition for identification holds for stationary and nonstationary data and some sort of instrumental variable estimators will have to be employed to yield consistency. However, nonstationarity does raise issues of inference if the rank of cointegration or direction of nonstationarity is not known *a priori*. The usual test statistics may not be chi-square distributed because of the presence of unit roots distributions. Classical instrumental variable estimators have to be modified to ensure valid inference.

## 1. Introduction

Let $\{\underline{w}_t\}$ be a sequence of time series observations of random variables. Multivariate vector autoregressive model (VAR) has been suggested as a useful tool to summarize the information contained in the data and to generate predictions (e.g. Hsiao [21, 22], Sims [50]). These models treat all variables as joint dependent and treat $\underline{w}_t$ as a function of its past values, $\underline{w}_{t-j}$. On the other hand, Cowles Commission approach assumes each equation in the system describes a behavioral or technological relations. An essential element of the Cowles Commission approach is to decompose $\underline{w}_t$ into $G$ endogenous variables, $\underline{y}_t$, and $K$ exogenous variables, $\underline{x}_t, \underline{w}'_t = (\underline{y}'_t, \underline{x}'_t), G + K = m$. The value of endogenous variables $\underline{y}_t$ are determined by the simultaneous interaction of the behavioral, technological or institutional relations in the model given the value of the exogenous variables, $\underline{x}_t$, and shock of the system (say, $\underline{\epsilon}_t$). The value of $\underline{x}_t$ is assumed to be determined by the forces outside of the model (e.g. Koopmans and Hood [19]). The Cowles Commission structural equation approach is also referred as a structural equations model (SEM). It has wide applications in education, psychology and econometrics, etc. (e.g. Browne and Arminger [6], Hood and Koopmans [19], Muthen [39, 40], Yuan and Bentler [59]). In this paper we will only focus on the aspects related to the time series analysis of a SEM.

Since the observed data can only provide information on conditional distribution of $\underline{y}_t$ given past values of $\underline{y}_{t-j}$ and current and past values of $\underline{x}_{t-j}$, there is an issue of if it is possible to infer from the data the true data generating process for the SEMs, which is referred to as an *identification* issue. Another issue for the SEMs is because of the joint dependency of $\underline{y}_t$, the regressors of an equation are correlated with the error (shock) of an equation which violates the condition for the regression method to be consistent. This is referred to as *simultaneity bias* issue. The theory and statistical properties of SEMs are well developed for stationary data (e.g. Amemiya [2], Intriligator, Boskin and Hsiao [30]).

[1]Department of Economics, University of Southern California, 3620 S. Vermont Ave. KAP300, Los Angeles, CA 90089, e-mail: chsiao@usc.edu

Nelson and Plosser [41] have shown that many economic and financial data contain unit roots, namely, most are integrated of order 1 or 2, I(1) or I(2). Theories for the time series analysis with unit roots have been derived by Anderson [4], Chan and Wei [7], Johansen [31, 32], Phillips [45], Phillips and Durlauf [46], Sims, Stock and Watson [51], Tiao and Tsay [57], etc. Among the major findings are that (i) $\underline{w}_t$ may be cointegrated in the sense that a linear combination of $I(d)$ variables may be of order $I(d-c)$, where $d$ and $c$ are positive numbers, say 1 (Granger and Weiss [14], Engle and Granger [11], Tiao and Box [54]); (ii) "Since these models (VAR) don't dichotomize variables into "endogenous" and "exogenous," the exclusion restrictions used to identify traditional simultaneous equations models make little sense" (Watson [58]); (iii) Time series regressions with integrated variables can behave very differently from those with stationary variables. Some of the estimated coefficients converge to their true values at the speed of $\sqrt{T}$ and are asymptotically normally distributed. Some converge to the true values at the speed of $T$ but have non-normal asymptotic distribution, and are asymptotically biased. Hence the Wald test statistics under the null may not be approximated by chi-square distributions (Chan and Wei [7], Sims, Stock and Watson [51], Tsay and Tiao [57]); (iv) Even though the $I(1)$ regressors may be correlated with the errors, the least squares regression consistently estimates the cointegrating relation, hence the simultaneity bias issues may be ignored (Phillips and Durlauf [46], Stock [52]).

In this paper we hope to review the recent advances in nonstationary time series analysis from the perspective of Cowles Commission Structural equation approach. In section 2 we discuss the relationships between a vector autoregressive model (VAR), a structural vector autoregressive model (SVAR), and Cowles Commission structural equations model (SEM). Section 3 discusses issues of estimating VAR with integrated variables. Section 4 discusses the least squares and instrumental variable estimators, in particular, the two stage least squares estimator (2SLS) for a SVAR. Section 5 discusses the modified and lag order augmented 2SLS estimators for SVAR. Conclusions are in Section 6.

## 2. Vector autoregression, structural vector autoregression and structural equations model

For ease of exposition, we shall assume that all elements of $\underline{w}_t$ are I(1) processes. We assume that $\underline{w}_t$ are generated by the following $p$-th order structural vector autoregressive process without intercept terms:[1]

$$(2.1) \qquad\qquad A(L)\underline{w}_t = \underline{\epsilon}_t$$

where $A(L) = A_0 + A_1 L + A_2 L^2 + \cdots + A_p L^p$. We assume that initial observations $\underline{w}_0, \underline{w}_{-1}, \ldots, \underline{w}_{-p}$ are available and

A.1: $A_0$ is nonsingular and $A_0 \neq I_m$, where $I_m$ denotes an $m$ rowed identity matrix.

---

[1]The introduciton of intercept terms complicates algebraic manipulation without changing the basic message. For detail, see [28].

A.2: The roots of $|A(L)| = 0$ are either 1 or outside the unit circle.

A.3: The $m \times 1$ error or innovation vector $\underline{\epsilon}_t$ is independently, identically distributed (i.i.d.) with mean zero, nonsingular covariance matrix $\Sigma_{\epsilon\epsilon}$ and finite fourth cumulants.

Premultiplying $A_0^{-1}$ to (2.1) yields the conventional VAR model of Johansen [31, 32], Phillips [45], Sims [50], Sims, Stock and Watson [51], Tsay and Tiao [57], etc.,

$$(2.2) \qquad \underline{w}_t = \Pi_1 \underline{w}_{t-1} + \cdots + \Pi_p \underline{w}_{t-p} + \underline{v}_t,$$

where $\Pi_j = -A_0^{-1} A_j, j = 1, \ldots, p$, and $\underline{v}_t = A_0^{-1} \epsilon_t$. The difference between (2.1) and (2.2) is that each equation in the former is supposed to describe a behavioral or technological relation while the latter is a *reduced form* relation. Eq. (2.2) is useful for generating prediction, but cannot be used for structural or policy analysis. For instance, $w_{1t}, w_{2t}, w_{3t}, w_{4t}$ may denote the price and quantity of a product, per capita income and raw material price, respectively. The first and second equations describe a demand relation which has quantity inversely related to price and positively related to income, and a supply relation which has price positively related to quantity and raw material price, respectively. Only (2.1) can provide information on demand and supply price elasticities but not (2.2). Equation (2.2) can only yield expected value of price and quantity given past $\underline{w}_{t-j}$.

Let $A = [A_0, A_1, \ldots, A_p]$ and define a $(p+1)m$-dimensional nonsingular matrix $M$ as

$$(2.3) \qquad M = \begin{bmatrix} I_m & I_m & \ldots & I_m \\ 0 & I_m & \ldots & I_m \\ 0 & 0 & \ldots & I_m \\ \cdot & \cdot & \cdots & \\ 0 & \cdot & 0 & I_m \end{bmatrix}.$$

Postmultiplying $A$ by $M$ yields an error-correction representation of (2.1),

$$(2.4) \qquad \sum_{j=0}^{p-1} A_j^* \bigtriangledown \underline{w}_{t-j} + A_p^* \underline{w}_{t-p} = \underline{\epsilon}_t,$$

where $\bigtriangledown = (1-L), A_j^* = \sum_{\ell=0}^{j} A_\ell, j = 0, 1, \ldots, p$. Let $A^* = [A_1^*, \ldots, A_p^*] = [\tilde{A}_1^*, A_p^*]$, then $A^* = AM$. The coefficient matrices $\tilde{A}_1^*$ and $A_p^*$ provide the implied *short-run* dynamics and *long-run* relations of the system (2.1) as defined in [26].[2]

Similarly, we can post-multiply (2.2) by $M$ to yield an *error-correction* representation of the reduced form (2.2)

$$(2.5) \qquad \bigtriangledown \underline{w}_t = \Pi_1^* \bigtriangledown \underline{w}_{t-1} + \cdots + \Pi_{p-1}^* \bigtriangledown \underline{w}_{t-p+1} + \Pi_p^* \underline{w}_{t-p} + \underline{v}_t,$$

where $\Pi_j = \sum_{i=1}^{j} \Pi_i - I_m$.

In this paper we are concerned with statistical inference of (2.1). If the roots of $|A(L)| = 0$ are all outside the unit circle, $\underline{w}_t$ is stationary. It is well known that the least squares estimator (LS) is inconsistent. The 2SLS and 3SLS using lagged $\underline{w}_t$ as instruments are consistent and asymptotically normally distributed (e.g. Amemiya

---

[2] The long-run and short-run dichotomization defined here is derived from (2.1). They are different from the those implied by Granger and Lin [13], Johansen [31, 32] or Pesaran, Shin and Smith [43], etc.

[2], Malinvaud [38]). Therefore, we shall assume that at least one root of $|A(L)| = 0$ is equal to 1. More specifically,[3]

A4:(a) $A_p^* = \underset{\sim}{\alpha}\beta'$ (or $\Pi_p^* = \underset{\sim}{\alpha}^*\beta^{*'}$) where $\underset{\sim}{\alpha}$ and $\beta$ (or $\underset{\sim}{\alpha}^*$ and $\beta^*$) are $m \times r$ matrices of full column rank $r, 0 \leq r \leq m-1$

   (b) $\underset{\sim}{\alpha}'_\perp J \underset{\sim}{\beta}_\perp$ or $(\underset{\sim}{\alpha}^{*'}_\perp J^* \beta_\perp^*)$ is nonsingular, where $J = \sum_{j=0}^{p-1} A_j^*$, (or $J^* = \sum_{j=0}^{p-1} \Pi_j^*$), $\underset{\sim}{\alpha}_\perp$ and $\underset{\sim}{\beta}_\perp$ (or $\underset{\sim}{\alpha}^*_\perp$ and $\beta^*_\perp$) are $m \times (m-r)$ matrices of full column rank such that $\underset{\sim}{\alpha}'_\perp \underset{\sim}{\alpha} = \underset{\sim}{0} = \beta'_\perp \beta$, (or $\underset{\sim}{\alpha}^{*'}_\perp \underset{\sim}{\alpha}^* = \underset{\sim}{0} = \beta^{*'}_\perp \beta$) (If $r = 0$, then we take $\underset{\sim}{\alpha}_\perp = I_m = \underset{\sim}{\beta}_\perp$.)

Under A1-A4, $\underset{\sim}{w}_t$ has $r$ *cointegrating* vectors (the columns of $\beta$) and $m - r$ unit roots. As shown by Johansen [31, 32] and Toda and Phillips [56] that A4 ensures that the *Granger representation theorem* (Engle and Granger [11]) applies, so that $\bigtriangledown \underset{\sim}{w}_t$ is stationary, $\beta' \underset{\sim}{w}_t$ is stationary, and $\underset{\sim}{w}_t$ is an I(1) process when $r < m$.

The cointegrating vectors $\beta$ provide information on the "long-run" or "equilibrium" state in which a dynamic system tends to converge over time after any of the variables in the system being perturbed by a shock, $\underset{\sim}{\alpha}$ transmits the deviation from such long-run relations, $\underset{\sim}{e}_t = \beta' \underset{\sim}{w}_t$, into each of $\underset{\sim}{w}_t$, and $\tilde{A}_1^*$ provides information on how soon such "equilibrium" is restored. In economics, the existence of long-run relationships and strength of attraction to such a state depends on the actions of a market or on government intervention. In this sense, the concept of *cointegration* has been applied in a variety of economic models including the relationships between capital and output; real wages and labor productivity; nominal exchange rate and relative prices, consumption and disposable income, long- and short-term interest rates, money velocity and interest rates, price of shares and dividends, production and sales, etc. (e.g. Banerjee, Dolado, Galbraith and Hendry [5], Hsiao, Shen and Fujiki [29], King, Plosser, Stock and Watson [33]).

Since the data only provide information of the conditional density of $\underset{\sim}{w}_t$ given past values of $\underset{\sim}{w}_{t-j}, j = 1, \ldots,$ there is an issue of if it is possible to derive (2.1) from (2.2) (or (2.4) from (2.5)). Without prior restrictions, there can be infinitely many different SVAR that yield identical (2.2). To see this we note that premultiplying (2.1) by any nonsingular constant matrix $F$ yields

$$(2.6) \qquad \tilde{A}_0 \underset{\sim}{w}_t + \tilde{A}_1 \underset{\sim}{w}_{t-1} + \cdots + \tilde{A}_p \underset{\sim}{w}_{t-p} = \bar{\underset{\sim}{\epsilon}}_t,$$

where $\tilde{A}_j = FA_j, \tilde{\underset{\sim}{\epsilon}}_t = F\underset{\sim}{\epsilon}_t$. Equations (2.1) and (2.5) yield identical (2.2) since $\tilde{A}_0^{-1}\tilde{A}_j = A_0^{-1}F^{-1}FA_j = \Pi_j, \underset{\sim}{v}_t = \tilde{A}_0^{-1}\tilde{\underset{\sim}{\epsilon}}_t = A_0^{-1}F^{-1}F\underset{\sim}{\epsilon}_t = A_0^{-1}\underset{\sim}{\epsilon}_t$. In other words, (2.1) and (2.5) are *observationally equivalent*.

An equation in (2.1) is identified if and only if the $g$-th row of *admissible transformation* matrix $F = (\underset{\sim}{f}'_g)$ takes the form that apart from the $g$th element being a nonzero constant, the rest are all zeros, i.e., $\underset{\sim}{f}'_g = (0, \ldots, 0, f_{gg}, 0, \ldots, 0)$ (e.g. Hsiao [23]). The transformation matrix $F$ is *admissible* if and only if (2.1) and (2.6) satisfy the same prior restrictions. Suppose that the $g$-th equation of (2.1) satisfies the prior restrictions $\underset{\sim}{a}'_g \Phi_g = \underset{\sim}{0}'$, where $\underset{\sim}{a}'_g$ denotes the $g$-th row of $A$ and $\Phi_g$ denotes a $(p+1)m \times R_g$ matrix with known elements. Let $\Phi_g^* = M^{-1}\Phi_g$, the existence of prior restrictions $\underset{\sim}{a}'_g \Phi_g = \underset{\sim}{0}'$ is equivalent to the existence of prior restrictions $\underset{\sim}{a}_g^{*'} \Phi_g^* = \underset{\sim}{0}'$, where $\underset{\sim}{a}_g^{*'}$ is the $g$-th row of $A^*$. It is shown by Hsiao [26] that

---

[3]Since $\Pi_p^* = A_0^{-1}A_p^*$, A4 implies that (a) $\Pi_p^* = \underset{\sim}{\alpha}^*\beta^{*'}$, where $\underset{\sim}{\alpha}^*$ and $\beta^*$ are $m \times r$ matrices of full column rank $r, 0 \leq r \leq m-1$, and (b) $\underset{\sim}{\alpha}_\perp^{*'} J^* \underset{\sim}{\beta}_\perp^*$ is nonsingular, where $J^* = \sum_{j=0}^{p-1} \Pi_j^*$.

**Theorem 2.1.** *Suppose that the g-th equation of (2.1) is subject to the prior restrictions $\underline{a}'_g \Phi_g = \underline{0}'$. A necessary and sufficient condition for the identification of the g-th equation of (2.1) or (2.4) is that*

$$(2.7) \qquad \operatorname{rank}(A\Phi_g) = m - 1,$$

*or*

$$(2.8) \qquad \operatorname{rank}(A^*\Phi_g^*) = m - 1.$$

Let $\underline{w}'_t = (\underline{y}'_t, \underline{x}'_t)$, where $\underline{y}'_t$ and $\underline{x}'_t$ are $1 \times G$ and $1 \times K$, respectively, and $G + K = m$. Let

$$A(L) = \begin{bmatrix} A_{11}(L) & A_{12}(L) \\ A_{21}(L) & A_{22}(L) \end{bmatrix},$$

and $\underline{\epsilon}'_t = (\underline{\epsilon}'_{1t}, \underline{\epsilon}'_{2t})$ be the conformable partitions. Cowles Commission decomposition of $\underline{w}_t$ into joint dependent variable variables $\underline{y}_t$ and exogenous variables $\underline{x}_t$ is equivalent to imposing the prior restrictions (Zellner and Palm [60]),

$$(2.9) \qquad A_{21}(L) \equiv \underline{0} \quad \text{and} \quad E\underline{\epsilon}_{1t}\underline{\epsilon}'_{2t} = \underline{0}.$$

The prior restrictions (2.9) restrict the *admissible transformation matrix $F$* to be block diagonal (e.g. Hsiao [23]). Therefore,

**Corollary 2.1.** *Under (2.9) and $\underline{a}'_g \Phi_g = \underline{0}'$, a necessary and sufficient condition for the identification of the g-th equation for $g \leq G$ is*

$$(2.10) \qquad \operatorname{rank}[(A_{11} \quad A_{12})\Phi_g] = G - 1,$$

*where $A_{11}$ and $A_{12}$ are conformable partitions of $A$.*

The identification condition (2.7) or (2.8) does not require any prior knowledge of the direction of nonstationarity or the rank of cointegration. As a matter of fact many macroeconometric models are identified without any prior knowledge of location of unit roots or rank of cointegration, (e.g. the Klein [34] interwar model and the large scale Wharton quarterly model (Klein and Evans [35]). Of course, if such information is available, it can improve the efficiency of system estimators and simplify the issues of inference considerably (e.g. King, Plosser, Stock and Watson [33]).

## 3. Inference in VAR (or reduced form)

Consider the $g$-th equation of (2.2),

$$(3.1) \qquad \underline{w}_g = X\underline{\pi}_g + \underline{v}_g,$$

where $\underline{w}_g$ is the $T \times 1$ vector of the $g$-th element of $\underline{w}_t$, $w_{gt}$, $X = (W_{-1}, \ldots, W_{-p})$, is the $T \times mp$ vector of $\underline{w}_{t-1}, \ldots, \underline{w}_{t-p}, \underline{\pi}_g$ is the corresponding vector of coefficients, and $\underline{v}_g$ is the $T \times 1$ vector of the $g$-th element of $\underline{v}_t$, $v_{gt}$.

Rewrite (3.1) in terms of linearly independent $I(0)$ and full rank $I(1)$ regressors $X_1^*$ and $X_2^*$, respectively, by postmultiplying a nonsingular transformation matrix

$M_x$ to $X$,[4] we have

$$
\begin{aligned}
(3.2) \qquad \underline{w}_g &= X M_x M_x^{-1} \underline{\pi}_g + \underline{v}_g \\
&= X^* \underline{\pi}_g^* + \underline{v}_g \\
&= (X_1^*, X_2^*) \begin{pmatrix} \underline{\pi}_{g1}^* \\ \underline{\pi}_{g2}^* \end{pmatrix} + \underline{v}_g,
\end{aligned}
$$

where $\underline{\pi}_g^* = M_x^{-1} \underline{\pi}_g = (\underline{\pi}_{g1}^{*'}, \underline{\pi}_{g2}^{*'})'$. The least squares estimator of (3.1) is equal to $M_x$ times the least squares estimator of (3.2),

$$
\begin{aligned}
(3.3) \qquad \hat{\underline{\pi}} &= (X'X)^{-1}(X'\underline{w}_g) \\
&= M_x (X^{*'} X^*)^{-1} X^{*'} \underline{w}_g \\
&= M_x [\underline{\pi}_g^* + (X^{*'} X^*)^{-1} X^{*'} \underline{v}_g].
\end{aligned}
$$

The statistical properties of (3.3) can be derived by making use of the fundamental functional central limit theorems proved by Chan and Wei [7], Phillips and Durlauf [46], etc.:

**Theorem 3.1.** *Let $\underline{\eta}_t$ be an $m \times 1$ vector of random variables with $E(\underline{\eta}_t \mid \underline{\eta}_{t-1}, \dots,) = \underline{0}$, $E(\underline{\eta}_t \underline{\eta}_t' \mid \underline{\eta}_{t-1}, \dots,) = I_m$, and bounded fourth moments. Let $F(L) = \sum_{j=0}^{\infty} F_j L^j$ and $G(L) = \sum_{j=0}^{\infty} G_j L^j$ with $\sum_{j=0}^{\infty} j \mid F_j \mid < \infty$ and $\sum_{j=0}^{\infty} j \mid G_j \mid < \infty$. Let $\xi_t = \sum_{s=1}^{t} \eta_s$, and let $B(r)$ denote an $m \times 1$ dimensional Brownian motion process.*

Then

(a) $T^{-1/2} \sum_{t=1}^{T} F(L) \underline{\eta}_t \Longrightarrow N(0, F(1)F(1)')$,
(b) $T^{-1} \sum_{t=1}^{T} \xi_{t-1} \eta_t \Longrightarrow \int B(r) dB(r)'$,
(c) $T^{-1} \sum_{t=1}^{T} \xi_t [F(L) \underline{\eta}_t]' \Longrightarrow F(1)' + \int B(r) dB(r)' F(1)'$,
(d) $T^{-1} \sum_{t=1}^{T} [F(L) \underline{\eta}_t][G(L) \underline{\eta}_t]' \longrightarrow \sum_{j=0}^{\infty} F_j G_j'$,
(e) $T^{-2} \sum_{t=1}^{T} \xi_t \xi_t' \Longrightarrow \int B(r) B(r)' dr$,

where to simplify notation $\int_0^1$ is denoted by $\int$ and $\longrightarrow$ and $\Longrightarrow$ denote convergence in probability and distribution of the associated probability measure, respectively.

Making use of theorem 3.1, it follows that

**Theorem 3.2.** *Under Assumptions A.1 - A.4, as $T \longrightarrow \infty$,*

$$
(3.4) \qquad \sqrt{T}(\hat{\underline{\pi}}_{g1}^* - \underline{\pi}_{g1}^*) \Longrightarrow N(\underline{0}, \sigma_{v_g}^2 M_{x_1 x_1}^*),
$$

$$
\begin{aligned}
(3.5) \qquad T(\pi_{g2}^* - \pi_{g2}^*) \Longrightarrow &\left( \int B_{x_2^*}(r) B_{x_2^*}(r)' dr \right)^{-1} \\
&\left( \int B_{x_2^*}(r) dB_{v_g}(r) \right).
\end{aligned}
$$

where $M_{x_1 x_1}^* = plim \frac{1}{T} \sum_{t=1}^{T} \underline{x}_{1t}^* \underline{x}_{1t}^{*'}$. Moreover, (3.4) and (3.5) are asymptotically independent.

---

[4] Such a transformation always exist. However, it does not need to be known *a priori*. The use of (3.2) is to facilitate the derivation of statistical properties of the estimators of (3.1) or (2.1).

The least squares estimator (3.3) is a linear combination of $\hat{\underline{\pi}}_{g1}^*$ and $\hat{\underline{\pi}}_{g2}^*$. Its limiting distribution is determined by the limiting distribution of the slower rate of $\hat{\underline{\pi}}_g^*$ included. Since the limiting distribution of $\hat{\underline{\pi}}_{g2}^*$ is nonstandard and involves a matrix unit distribution, the usual Wald test statistic under the null may not be approximated by the chi-square distribution if the null hypothesis involves coefficients in the direction of nonstationarity (e.g. Dolado and Lutkepohl [9], Sims, Stock and Watson [51], Tsay and Tiao [57]). On the other hand, if $\underline{w}_t$ is cointegrated and the rank of cointegration is known *a priori*, Ahn and Reinsel [1] and Johansen [31, 32] using the reduced rank framework proposed by Anderson [3] have shown that the coefficients of cointegration vectors are asymptotically mixed normal, hence there will be no inference problem. The Wald test statistics constructed from the reduced rank regression will again be asymptotically chi-square distributed. This is because imposing the reduced rank condition is equivalent to avoid estimating the unit roots in the system.

Unfortunately, as discussed in section 2, prior information on the rank of cointegration or direction of nonstationarity is usually lacking. One way to deal with it is to pretest the data for the presence of cointegration and the rank of cointegration, then apply the reduced rank regression of Ahn and Reinsel [1] or Johansen [31, 32]. However, statistic tests for the rank of cointegration have very poor finite sample performance (e.g. Stock [53]). The first stage unit root test and second stage cointegration test can induce substantial size distortion. For instance, Elliott and Stock [10] consider a bivariate problem in which there is uncertainty about whether the regressor has a unit root. In their Monte Carlo simulation they find that unit root pretests can induce substantial size distortions in the second-stage test. If the innovations of the regressors and the second-stage regression error are correlated, the first-stage Dickey-Fuller [8] *t*-statistic and the second-stage *t*-statistic will be dependent so the size of the second stage in this two-stage procedure cannot be controlled, even asymptotically. Many other Monte Carlo studies also show that serious size and power distortions arise and the number of linearly independent cointegrating vectors tend to be overestimated as the dimension of the system increases relative to the time dimension (e.g. Ho and Sorensen [18], Gonzalo and Pitarakis [12]).

Another way is to correct the miscentering and skewness of the limiting distribution of the least squares estimator due to the "endogeneities" of the predetermined integrated regressors (e.g. Park [42], Phillips [44], Phillips and Hansen [47], Robinson and Hualde [49]). However, since the rank of cointegration and direction of nonstationarity are unknown, Phillips [45] proposes to deal with potential endogeneities by making a correction of the least squares regression formula that adjusts for whatever endogeneities there may be in the predetermined variables that is due to their nonstationarity by transforming the dependent variables $\underline{w}_t$ into

$$(3.6) \qquad \underline{w}_t^+ = \underline{w}_t - \Omega_{v\triangledown w} \Omega_{\triangledown w\triangledown w}^- \triangledown \underline{w}_t,$$

where $\Omega_{\triangledown w\triangledown w} = \sum_{j=-\infty}^{\infty} E(\triangledown \underline{w}_t \triangledown \underline{w}_{t-j}')$, $\Omega_{v\triangledown w} = \sum_{j=-\infty}^{\infty} E(\underline{v}_t \triangledown \underline{w}_{t-j}')$ and $\Omega_{\triangledown w\triangledown w}^-$ denotes the Moore-Penrose generalized inverse.[5] Using $\underline{w}_t^+$ in place of $\underline{w}_t$ in (2.2) is equivalent to modifying the error term from $\underline{v}_t$ to $\underline{v}_t - \Omega_{\triangledown w}\Omega_{\triangledown w\triangledown w}^- \triangledown \underline{w}_t$, which now becomes serially correlated because $\triangledown \underline{w}_t$ is serially correlated. To correct for this order $(1/T)$ serial correlation bias term, Phillips [45] suggests further adding $(X'X)^{-1}(0, T\Delta_{v\triangledown w}^{+'})$ to the least squares regression estimator of $\underline{w}_t^+$ on $\triangledown \underline{w}_{t-1}, \ldots, \triangledown \underline{w}_{t-p+1}, \underline{w}_{t-p}$, where $\Delta_{v\triangledown w}^+ = \Omega_{v\triangledown w}\Omega_{\triangledown w\triangledown w}^- \Delta_{\triangledown w\triangledown w}$, and $\Delta_{uv}$ denotes the one-sided long-run covariances of two sets of $I(0)$ variables $(u_t, v_t)$, $\Delta_{uv} = \sum_{j=0}^{\infty} \Gamma_{uv}(j)$ where $\Gamma_{uv}(j) = E\underline{u}_t\underline{v}_{t-j}$.[6] Consistent estimates of $\Omega_{uv}$ or $\Delta_{uv}$ can be obtained by using Kernel method (e.g. Hannan [15], Priestley [48]).

---

[5]If $\underline{w}_t$ are cointegrated, $\Omega_{\triangledown w\triangledown w}$ does not have full rank.

$$(3.7) \qquad \hat{\Omega}_{uv} = \sum_{j=-T+1}^{T-1} h(j/K)\hat{\Gamma}_{uv}(j),$$

$$(3.8) \qquad \hat{\Delta}_{uv} = \sum_{j=0}^{T-1} h(j/K)\hat{\Gamma}_{uv}(j),$$

where $\hat{\Gamma}_{uv}(j)$ is a consistent sample covariance estimator of $\Gamma_{uv}(j)$, and $h(\cdot)$ is a kernel function and $K$ is a lag truncation or bandwidth parameter. Assuming that

**Assumption 3.1.** The kernel function $h(\cdot) : R \longrightarrow [-1, 1]$ is a twice continuously differentiable even function with:

   (a) $h(0) = 1, h'(0) = 0, h''(0) \neq 0$; and either
   (b) $h(x) = 0, \mid x \mid \geq 1$, with $\lim_{|x| \longrightarrow 1} \frac{h(x)}{(1-|x|)^2} =$ constant, or
   (b') $h(x) = O((1-x)^2)$, as $\mid x \mid \longrightarrow 1$.

**Assumption 3.2.** The bandwidth parameter $K$ in the kernel estimates (3.7) and (3.8) has an expansion rate $K \sim c_T T^k$ for some $k \in (1/4, 2/3)$ and for some slowly varying function $c_T$ and thus $K/T^{2/3} + T^{1/4}/K \longrightarrow 0$ and $K^{4/T} \longrightarrow \infty$ as $T \longrightarrow \infty$.

Phillips [45] shows that the modified least squares estimates are either asymptotically normally distributed or mixed normal. However, because the direction of nonstationarity is unknown, the conditional covariance matrix cannot be derived. Therefore, if the test statistic involves some of the coefficients of nonstationary variables, the limiting distribution becomes a mixture of chi-squares variates with the weights between 0 and 1. In other words, if tests based on chi-square distribution rejects the null with significance level $\alpha$, then the test rejects the null with significance level less than $\alpha$. In other words, tests based on chi-square distribution provides a conservative test.

Toda and Yamamoto [55] have suggested a lag-order augmented approach to circumscribe the issue of non-standard distributions associated with integrated regressors by overfitting a VAR with additional $d_{\max}$ lags where $d_{\max}$ denotes the maximum order of integration suspected. In our case, $d_{\max} = 1$. In other words, instead of estimating (2.2), we estimate

$$(3.9) \qquad \underline{w}_t = \Pi_1\underline{w}_{t-1} + \cdots + \Pi_p\underline{w}_{t-p} + \Pi_{p+1}\underline{w}_{t-p-1} + \underline{v}_t,$$

Since we know *a priori*, $\Pi_{p+1} \equiv 0$, we are only interested in the estimates of $\Pi_j, j = 1, \ldots, p$. The limiting distributions of the least squares estimates of (3.9) can be derived from the limiting distributions of the least squares estimates of (the error-correction form),

$$(3.10) \qquad \underline{w}_t = \Pi_1^* \bigtriangledown \underline{w}_{t-1} + \cdots + \Pi_p^* \bigtriangledown \underline{w}_{t-p} + \Pi_{p+1}^*\underline{w}_{t-p-1} + \underline{v}_t,$$

because $\Pi_j^* = \sum_{i=1}^{j} \Pi_i, j = 1, \ldots, p+1$ or $\Pi_j = \Pi_j^* - \Pi_{j-1}^*$ where $\Pi_0^* \equiv \underline{0}$. Since $\Pi_j^*, j = 1, \ldots, p$ are coefficients of stationary regressors, Theorem 3.2 shows that the least squares estimates of $\Pi_j^*$, $j = 1, \ldots, p$ converge to the true values at the speed of $\sqrt{T}$ and are asymptotically normally distributed. Only the least squares estimates of $\Pi_{p+1}^*$ may be $T$-convergent and have non-normal limiting distributions. However, since we know *a priori* that $\Pi_{p+1} = \underline{0}$, our interest is only in $\Pi_j$, $j = 1, \ldots, p$. The least squares regression of (3.9) yields $\hat{\Pi}_j = \hat{\Pi}_j^* - \hat{\Pi}_{j-1}^*$, $j = 1, \ldots, p$, therefore, they are asymptotically normally distributed. Wald test statistics of the null hypothesis constructed from regression estimates of (3.9) will again be asymptotically chi-square distributed.

---

[6] Under A.3, $\Delta_{v\bigtriangledown w} = \underline{0}$.

Phillips [45] modified estimator maintains the $T$-convergence part of the coefficients associated with full rank integrated regressors. The Toda-Yamamoto [55] lag order augmented estimator is only $\sqrt{T}$-convergent. So Phillips [45] modified estimator is likely to be asymptotically more efficient. However, computationally, the Phillips modified estimator is much more complicated than the lag order augmented estimator. Moreover, test statistics constructed from the modified estimators can only give the bounds of the size of the test because the conditional variance is unknown, while test statistics constructed from the lag order augmented estimator asymptotically yield the exact size.

## 4. Least squares and two stage least squares estimation of SVAR

For ease of exposition, we assume that prior information is in the form of excluding certain variables, both current and lagged, from an equation. Let the $g$-th equation of (2.1) be written as

$$(4.1) \qquad \underline{w}_g = Z_g \underline{\delta}_g + \underline{\epsilon}_g,$$

where $\underline{w}_g$ and $\underline{\epsilon}_g$ denote the $T \times 1$ vectors of $(w_{g1}, \ldots, w_{gT})'$ and $(\epsilon_{g1}, \ldots, \epsilon_{gT})'$, respectively, and $Z_g$ denotes the $T \times [(p+1)g_\Delta - 1]$ dimensional matrix of $g_\Delta$ included current and lagged variables of $\underline{w}_t$.

The least squares estimator of (4.1) is given by

$$(4.2) \qquad \hat{\underline{\delta}}_{g,\ell s} = (Z_g' Z_g)^{-1} Z_g' \underline{w}_g$$

Phillips and Durlauf [46] and Stock [52] have shown that the least squares estimator with integrated regressors is consistent even when the regressors and the errors are correlated. However, the basic assumption underlying their result is that the regressors are not cointegrated. In a dynamic framework even though $\underline{w}_{t-j}$ are $I(1)$, the current and lagged variables are trivially cointegrated. It was shown in [21] when contemporaneous joint dependent variables also appear as explanatory variables in (4.1), applying least squares method to (4.1) does not yield consistent estimator for $\underline{\delta}_g$. To see this, let $M_g$ be the nonsingular transformation matrix that transforms $Z_g$ into $Z_g^* = Z_g M_g = (Z_{g1}^*, Z_{g2}^*)$, where $Z_{g1}^*$ denotes the $\ell_g$-dimensional linearly independent $I(0)$ variables and $Z_{g2}^*$ denotes the $T$ observations of $b_g$ full rank $I(1)$ variables,[7] then

$$(4.3) \qquad \begin{aligned} \underline{w}_g &= Z_g M_g M_g^{-1} \underline{\delta}_g + \underline{\epsilon}_g \\ &= Z_g^* \underline{\delta}_g^* + \underline{\epsilon}_g \end{aligned}$$

where $\underline{\delta}_g^* = M_g^{-1} \underline{\delta}_g = (\underline{\delta}_{g1}^{*'}, \underline{\delta}_{g2}^{*'})'$ with $\underline{\delta}_{g1}^*$ and $\underline{\delta}_{g2}^*$ denoting the $\ell_g \times 1$ and $b_g \times 1$ vector, respectively. Such transformation always exists. For instance, if no cointegrating relation exists among the included $\underline{w}_t$, say $\tilde{w}_{gt}$, then $b_g$ equals the dimension of included joint dependent variables, $g_\Delta$, and $Z_{g2}^*$ consists of the first differenced current and $p-1$ lagged included variables, $Z_{g2}^*$ is simply the $T \times b_g$ (or $T \times g_\Delta$) included $\underline{\tilde{w}}_{gt}$ lagged by $p$ periods, $\underline{\tilde{w}}_{g,t-p}$. On the other hand, if there exists $g_\Delta - b_g$ linearly independent cointegrating relations among the $g_\Delta$ included variables, $\tilde{w}_{gt}$, then $Z_{g1}^*$ consists of the current and $p-1$ lagged $\bigtriangledown \tilde{w}_{gt}$ and $\tilde{W}_{g,-p}d_g$ cointegrating relations, where $\tilde{W}_{g,-p}$ is $T \times g_\Delta$ matrix of included $\underline{\tilde{w}}_{g,t-p}$, $d_g$ is $g_\Delta \times (g_\Delta - b_g)$ of constants, and $Z_{g2}^*$ consists of the $T$ observed $b_g$ full rank $I(1)$ variables $\tilde{W}_{g2,-p}$.

---

[7] By full rank I(1) variables we mean that there is no cointegrating relation among $Z_{g2}^*$.

The least squares estimator (4.2) can be written as $\hat{\underline{\delta}}_{g,\ell s} = M_g \hat{\underline{\delta}}_{g,\ell s}^*$, where $\hat{\underline{\delta}}_{g,\ell s}^*$ denotes the least squares estimator of (4.3). Using Theorem 3.1, one can show that $\frac{1}{T}Z_{g1}^{*'}Z_{g1}^* \longrightarrow M_{z_{g1}z_{g1}}^*$, $T^{-2/3}Z_{g1}^{*'}Z_{g2}^* \longrightarrow 0$, $\frac{1}{T^2}Z_{g2}^{*'}Z_{g2}^* \Longrightarrow M_{z_{g2}z_{g2}}^*$, $\frac{1}{T^2}Z_{g2}^{*'}\underline{\epsilon}_g \longrightarrow 0$, $\frac{1}{T}Z_{g1}^{*'}\underline{\epsilon}_g \longrightarrow \underline{b}$, where $\underline{b} = [E(\underline{\epsilon}_{gt}\underline{\tilde{w}}_{gt}'), \underline{0}']' = [(A_0^{-1}\sum_{\epsilon\epsilon,g})_g', \underline{0}']'$, $\sum_{\epsilon\epsilon}$ is the $g$-th column of $\sum_{\epsilon\epsilon}$ and $(A_0^{-1}\sum_{\epsilon\epsilon,g})_g$ is the $(g_\Delta - 1) \times 1$ subvector of $A_0^{-1}\sum_{\epsilon\epsilon,g}$ that corresponds to the $g_\Delta - 1$ included variables $\tilde{w}_{gt}$ in the $g$-th equation, and $M_{z_{g1}z_{g1}}^*$ and $M_{z_{g2}z_{g2}}^*$ are nonsingular. It follows that

$$\hat{\underline{\delta}}_{g,\ell s}^* = \begin{bmatrix} \hat{\underline{\delta}}_{g1,\ell s}^* \\ \hat{\underline{\delta}}_{g2,\ell s}^* \end{bmatrix} \longrightarrow \begin{bmatrix} \underline{\delta}_{g1}^* \\ \underline{\delta}_{g2}^* \end{bmatrix} + \begin{bmatrix} \underline{b} \\ \underline{0} \end{bmatrix}. \tag{4.4}$$

Although the coefficients of $Z_{g2}^*$ can be consistently estimated, the coefficients of $Z_{g1}^*$ cannot. Since $\hat{\underline{\delta}}_{g,\ell s}$ is a linear combination of $\hat{\underline{\delta}}_{g1,\ell s}^*$ and $\hat{\underline{\delta}}_{g2,\ell s}^*$, $\hat{\underline{\delta}}_{g,\ell s}$ is inconsistent.

When the errors and regressors are correlated, a standard procedure is to use instrumental variable method. Using lagged variables as instruments, the two stage least squares estimator of $\underline{\delta}_g$ is given by

$$\hat{\underline{\delta}}_{g,2SLS} = [Z_g'X(X'X)^{-1}X'Z_g]^{-1}[Z_g'X(X'X)^{-1}Z_g'\underline{w}_g], \tag{4.5}$$

where $X = (W_{-1}, W_{-2}, \ldots, W_{-p})$ and $W_{-j}$ denotes the $T \times m$ matrix representation of $\underline{w}_{t-j}$. Transforming $X$ into linearly independent $I(0)$ and full rank $I(1)$ processes, $X_1^*$ and $X_2^*$, respectively, by $M_x$, $XM_x = [X_1^*, X_2^*]$, the 2SLS estimator (4.5) is equal to $M_g\hat{\underline{\delta}}_{g,2SLS}^*$, where

$$\hat{\underline{\delta}}_{g,2SLS}^* = [Z_g^{*'}X^*(X^{*'}X^*)^{-1}X^{*'}Z_g^*]^{-1}[Z_g^{*'}X^*(X^{*'}X^*)^{-1}X^{*'}\underline{w}_g] \tag{4.6}$$

Since $\frac{1}{T^2}Z_{g1}^{*'}X_2^* \longrightarrow 0$, $\frac{1}{T}Z_{g2}^{*'}X_1^* \Longrightarrow M_{z_{g2}x_1}^*$, $\frac{1}{T^2}Z_{g2}^{*'}X_1^* \longrightarrow 0$, $\frac{1}{T}X_1^{*'}X_1^* \longrightarrow M_{x_1x_1}^*$, $\frac{1}{T}X_1^{*'}X_2^* \Longrightarrow M_{x_1x_2}^*$ $\frac{1}{T^2}X_1^{*'}X_2^* \longrightarrow 0$, $\frac{1}{T^2}X_2^{*'}X_2^* \Longrightarrow M_{x_2x_2}^*$ $\frac{1}{T}X_1^{*'}\underline{\epsilon}_g \longrightarrow 0$, and $\frac{1}{T^2}X_2^{*'}\underline{\epsilon}_g \longrightarrow 0$, and $M_{x_2x_2}^*$ are nonsingular, it follows that $\hat{\underline{\delta}}_{g,2SLS}^*$ converges to $\underline{\delta}_g^*$. Hence the 2SLS estimator of $\underline{\delta}_g$ is consistent.

Let $H_g = \begin{bmatrix} T^{-\frac{1}{2}}I_{\ell_g} & 0 \\ 0 & T^{-1}I_{b_g} \end{bmatrix}$ and $H_x = \begin{bmatrix} T^{-\frac{1}{2}}I_{\ell^*} & 0 \\ 0 & T^{-1}I_{b^*} \end{bmatrix}$, where $\ell^*$ and $b^*$ are the column dimensions of $X_1^*$ and $X_2^*$ respectively. Under assumptions A.1 - A.4, as $T \longrightarrow \infty$,

$$\begin{aligned} H_g^{-1}(\hat{\underline{\delta}}_{g,2SLS}^* - \underline{\delta}_g^*) &= \begin{bmatrix} \sqrt{T}(\hat{\underline{\delta}}_{g1,2SLS}^* - \underline{\delta}_{g1}^*) \\ T(\hat{\underline{\delta}}_{g2,2SLS}^* - \underline{\delta}_{g2}^*) \end{bmatrix} \\ &\Longrightarrow \begin{bmatrix} (M_{z_{g1}x_1}^* M_{x_1x_1}^{*-1} M_{x_1z_{g1}}^*)^{-1}(M_{z_{g1}x_1}^* M_{x_1x_1}^{*-1} \cdot T^{-1/2}X_1^{*'}\underline{\epsilon}_g) \\ (M_{z_{g2}x_2}^* M_{x_2x_2}^{*-1} M_{x_2z_{g2}}^*)^{-1}(M_{z_{g2}x_2}^* M_{x_2x_2}^{*-1} \cdot T^{-1}X_2^{*'}\underline{\epsilon}_g) \end{bmatrix}. \end{aligned} \tag{4.7}$$

By theorem 3.1, we have

$$\frac{1}{\sqrt{T}}X_1^{*'}\underline{\epsilon}_g \Longrightarrow N(\underline{0}, \sigma_g^2 M_{x_1x_1}^*). \tag{4.8}$$

and

$$(4.9) \qquad \frac{1}{T} X_2^{*\prime} \epsilon_g \Longrightarrow \int B_{x_2^*} dB_{\epsilon_g},$$

where $B_{\epsilon_g}$ denotes the Brownian motion of $\epsilon_{gt}$ with variance $\sigma_g^2$, $B_{x_2^*}$ denotes a $b^* \times 1$ vector Brownian motion of $\bigtriangledown \underline{x}_{2t}^*$ with covariance matrix $\Omega_{\bigtriangledown x_2^* \bigtriangledown x_2^*}$ where $\Omega_{\bigtriangledown x_2^* \bigtriangledown x_2^*}$ is the long-run covariance matrix of $\bigtriangledown \underline{x}_{2t}^*$. The Brownian motion $B_{x_2^*}^*$ and $B_{\epsilon_g}$ are not independent because $\epsilon_{gt}$ and $\underline{v}_t$ are contemporaneously correlated. Following Phillips [44], we can decompose the right hand side of (4.9) into two terms as

$$(4.10) \qquad \int B_{x_2^*} dB_{\epsilon_g \cdot x_2^*} + \int B_{x_2^*} \Omega_{\epsilon_g \bigtriangledown x_2^*} \Omega_{\bigtriangledown x_2^* \bigtriangledown x_2^*}^{-1} dB_{x_2^*},$$

where $B_{\epsilon_g \cdot x_2^*} = B_{\epsilon_g} - \Omega_{\epsilon_g \bigtriangledown x_2^*} \Omega_{\bigtriangledown x_2^* \bigtriangledown x_2^*}^{-1} B_{x_2^*} \equiv BM(\sigma_{g,\bigtriangledown x_2^*}^2)$ with $\sigma_{g,\bigtriangledown x_2^*}^2 = \sigma_g^2 - \Omega_{\epsilon_g \bigtriangledown x_2^*} \Omega_{\bigtriangledown x_2^* \bigtriangledown x_2^*}^{-1} \Omega_{\bigtriangledown x_2^* \epsilon_g}$, and $\Omega_{\epsilon_g \bigtriangledown x_2^*}$ denotes the long-run covariance between $\epsilon_g$ and $\bigtriangledown \underline{x}_2^*$. The first term of (4.10) is a mixed normal. The second term involves a matrix unit root distribution that arises from using lagged $\underline{w}$ as instruments when $\underline{w}$ is I(1) and the contemporaneous correlation between $\epsilon_{gt}$ and $\underline{w}_t$ is nonzero. The "long-run endogeneity" of the nonstationary instruments $X_2^*$ leads to a skewness of the limiting distribution of $\hat{\underline{\delta}}_{g,2SLS}^*$ and its dependence on nuisance parameters that are impossible to eliminate by the 2SLS. Therefore,

**Theorem 4.1.** *Under A.1 - A.4 the 2SLS estimator of $\delta_g^*$ is consistent and*

$$(4.11) \qquad \sqrt{T}(\hat{\delta}_{g1,2SLS}^* - \delta_{g1}^*) \Longrightarrow N(\underline{0}, \sigma_g^2 (M_{z_{g1}x_1}^* M_{x_1x_1}^{*-1} M_{x_1z_{g1}}^*)^{-1}),$$

$$
\begin{aligned}
(4.12) \qquad T(\hat{\delta}_{g2,2SLS}^* - \delta_{g2}^*) \Longrightarrow & \left\{ \int B_{z_{g2}^*} B'_{x_2^*} dr (\int B_{x_2^*} B'_{x_2^*} dr)^{-1} \int B_{x_2^*} B'_{z_{g2}^*} dr \right\}^{-1} \\
& \left\{ \int B_{z_{g2}^*} B'_{x_2^*} dr (\int B_{x_2^*} B'_{x_2^*} dr)^{-1} \right. \\
& \times \left. \left[ \int B_{x_2^*} dB_{\epsilon_g \cdot x_2^*} + \int B_{x_2^*} \Omega_{\epsilon_g \bigtriangledown x_2^*} \Omega_{\bigtriangledown x_2^* \bigtriangledown x_2^*}^{-1} dB_{x_2^*} \right] \right\},
\end{aligned}
$$

*where $B_{z_{g2}^*}$ denotes a $b_g \times 1$ vector Brownian motion of $\bigtriangledown z_{g2,t}^*$ which appears in the $g$-th equation. The distributions of (4.11) and (4.12) are asymptotically independent.*

Theorem 4.1 suggests that inference about the null hypothesis $P\delta_g = \underline{c}$ can be tricky, where $P$ and $\underline{c}$ are known matrix and vector of proper dimensions. If $\sqrt{T}P(\hat{\underline{\delta}}_{g,2SLS} - \underline{\delta}_g)$ has a nonsingular covariance matrix, the limiting distribution of $P\hat{\underline{\delta}}_g$ is determined by the limiting distribution of $\hat{\underline{\delta}}_{g1}^*$, hence the Wald test statistic

$$(4.13) \qquad (\hat{\underline{\delta}}_{g,2SLS} - \underline{\delta}_g)' P' \text{ Cov } (P\hat{\underline{\delta}}_{g,2SLS})^{-1} P(\hat{\underline{\delta}}_{g,2SLS} - \underline{\delta}_g)$$

under the null will be asymptotically chi-square distributed. On the other hand, if $\sqrt{T}P(\hat{\underline{\delta}}_{g,2SLS} - \underline{\delta}_g)$ has a singular covariance matrix, it means that there exists a nonsingular matrix $L$ such that

$$(4.14) \qquad LP\delta_g = LP^* \underline{\delta}_g^* = \begin{bmatrix} \tilde{P}_{11} & \tilde{P}_{12} \\ \underline{0} & \tilde{P}_{22} \end{bmatrix} \begin{bmatrix} \underline{\delta}_{g1}^* \\ \underline{\delta}_{g2}^* \end{bmatrix}$$

with nonzero $\tilde{P}_{22}$. Then

$$(P\hat{\underline{\delta}}_{g,2SLS} - \underline{c})' \text{ Cov } (P\hat{\underline{\delta}}_{g,2SLS})^{-1}(P\hat{\underline{\delta}}_{g,2SLS} - \underline{c})$$

$$= \left\{ \begin{bmatrix} \tilde{P}_{11} & \tilde{P}_{12} \\ \underline{0} & \tilde{P}_{22} \end{bmatrix} \begin{bmatrix} \hat{\underline{\delta}}^*_{g1,2SLS} \\ \hat{\underline{\delta}}^*_{g2,2SLS} \end{bmatrix} - L\underline{c} \right\}' \text{ Cov } (LP\hat{\underline{\delta}}_{g,2SLS})^{-1}$$

$$\times \left\{ \begin{bmatrix} \tilde{P}_{11} & \tilde{P}_{12} \\ \underline{0} & \tilde{P}_{22} \end{bmatrix} \begin{bmatrix} \hat{\underline{\delta}}^*_{g1,2SLS} \\ \hat{\underline{\delta}}^*_{g2,2SLS} \end{bmatrix} - L\underline{c} \right\}$$

(4.15)
$$\Longrightarrow T(\tilde{P}_{11}\hat{\underline{\delta}}^*_{g1,2SLS} + \tilde{P}_{12}\hat{\underline{\delta}}^*_{g2,2SLS} - \tilde{\underline{c}}_1)' \text{ Cov } (\sqrt{T}\tilde{P}_{11}\hat{\underline{\delta}}^*_{g1,2SLS})^{-1}$$

$$\times (\tilde{P}_{11}\hat{\underline{\delta}}^*_{g1,2SLS} + \tilde{P}_{12}\hat{\underline{\delta}}^*_{g2,2SLS} - \tilde{\underline{c}}_1)$$

$$+ T^2(\tilde{P}_{22}\hat{\underline{\delta}}^*_{g2,2SLS} - \tilde{\underline{c}}_2)' \text{ Cov } (T\tilde{P}_{22}\hat{\underline{\delta}}^*_{g2,2SLS})^{-1}(\tilde{P}_{22}\hat{\underline{\delta}}^*_{g2,2SLS} - \tilde{\underline{c}}_2),$$

where $L\underline{c} = (\tilde{\underline{c}}'_1, \tilde{\underline{c}}'_2)'$. The first term on the right hand side of (4.15) is asymptotically chi-square distributed. The second term, according to Theorem 3.1 has a nonstandard distribution. Hence (4.15) is not asymptotically chi-square distributed.

If there exists prior information that satisfies (2.9) and $\underline{w}_1$ and $\underline{w}_2$ are cointegrated with $x^*_2$ contained in $\underline{w}_2$, it was shown by Hsiao [22] that the 2SLS converges to a mixed normal distribution. Then the Wald test statistic (4.13) can again be approximated by a chi-square distribution. When variables cannot be dichotomized into "endogenous" and "exogenous", if we do not know the direction of nonstationarity, nor the rank of cointegration, we will not be able to know *a priori* if $P_{22}$ is a zero matrix, hence if (4.13) may be approximated by a chi-square distribution.

## 5. Modified and lag order augmented 2SLS estimators

We note that just like the least squares estimator for the VAR model, the application of 2SLS does not provide asymptotically normal or mixed normal estimator because of the long-run endogeneities between lagged I(1) instruments and the (current) shocks of the system. But if we can condition on the innovations driving the common trends it will allow us to establish the independence between Brownian motion of the errors of the conditional system involving the cointegrating relations and the innovations driving the common trends. The idea of the modified 2SLS estimator is to apply the 2SLS method to the equation conditional on the innovations driving the common trends. Unfortunately, the direction of nonstationarity is generally unknown. Neither does the identification condition given by Theorem 2.1 requires such knowledge. In the event that such knowledge is unavailable, Hsiao and Wang [27] propose to generalize Phillips [45] fully modified VAR estimator to the 2SLS estimator.

Rewrite (4.1) as

(5.1)
$$\underline{w}_g = Z_g \tilde{M}_g \tilde{M}_g^{-1} \underline{\delta}_g + \underline{\epsilon}_g$$
$$= (Z^{**}_{g1} \quad Z^{**}_{g2}) \begin{pmatrix} \underline{\delta}^{**}_{g1} \\ \underline{\delta}^{**}_{g2} \end{pmatrix} + \underline{\epsilon}_g$$
$$= Z^{**}_g \underline{\delta}^{**}_g + \underline{\epsilon}_g$$

where $Z^{**}_g = Z_g \tilde{M}_g = (Z^{**}_{g1}, Z^{**}_{g2})$, $Z^{**}_{g1} = (\nabla W_g, \nabla \tilde{W}_{g,-1}, \ldots, \nabla \tilde{W}_{g,-p+1})$, $Z^{**}_{g2} = \tilde{W}_{g,-p}, \underline{\delta}^{**}_g = \tilde{M}_g^{-1}\underline{\delta}_g$, $\nabla \tilde{W}_{g,-j}$ denoting the $T \times g_\Delta$ stacked first difference of the included variable $\nabla \tilde{\underline{w}}_{g,t-j}$ and $\nabla W_g$ denoting the $T \times (g_\Delta - 1)$ first difference of the included variables $\triangledown \tilde{\underline{w}}_{gt}$ excluding $\triangledown w_{gt}$. The decomposition $(Z_{g1}^{**}, Z_{g2}^{**})$ and $\underline{\delta}_g^{**} = (\delta_{g1}^{**\prime}, \delta_{g2}^{**\prime})'$ are identical to $(Z_{g1}^*, Z_{g2}^*)$ if there is no cointegrating relations among $\tilde{\underline{w}}_{gt}$, $\underline{d}_g = \underline{0}$. Unlike $(Z_{g1}^*, Z_{g2}^*)$, $(Z_{g1}^{**}, Z_{g2}^{**})$ are well defined and observable. When $Z_{g1}^* \neq Z_{g1}^{**}$, there exists a nonsingular transformation matrix $D_g$ such that $(Z_{g1}^{**}, Z_{g2}^{**})D_g = (Z_{g1}^*, Z_{g2}^*)$. Then

$$(5.2) \qquad \underline{\delta}_g^* = D_g^{-1} \underline{\delta}_g^{**}.$$

Let

$$(5.3) \qquad C_g = (W_{-p}' \triangledown W_{-p} - T\Delta_{\triangledown w \triangledown w})\Omega_{\triangledown w \triangledown w}^{-} \Omega_{\triangledown w \epsilon_g},$$

where $\Omega_{uv}$ and $\Delta_{uv}$ denote the long-run covariance and the one-sided long-run covariance matrix of two sets of I(0) variables, $(\underline{u}_t, \underline{v}_t)$,

$$(5.4) \qquad \Omega_{uv} = \sum_{j=-\infty}^{\infty} \Gamma_{uv}(j),$$

and

$$(5.5) \qquad \Delta_{uv} = \sum_{j=0}^{\infty} \Gamma_{uv}(j),$$

where $\Gamma_{uv}(j) = E\underline{u}_t \underline{v}_{t-j}'$. Let

$$(5.6) \qquad \hat{C}_g = (W_{-p}' \triangledown W_{-p} - T\hat{\Delta}_{\triangledown w \triangledown w})\hat{\Omega}_{\triangledown w \triangledown w}^{-1} \hat{\Omega}_{\triangledown w \epsilon_g},$$

where $\hat{\Omega}_{uv}$ and $\hat{\Delta}_{uv}$ are the kernel estimates of $\Omega_{uv}$ and $\Delta_{uv}$, such as (3.7) and (3.8). A modified 2SLS estimator following Phillips [45] fully modified VAR estimator can be defined as

$$(5.7) \qquad \begin{aligned} \hat{\tilde{\delta}}_{g,m2SLS}^{**} = & \left\{ Z_g^{**\prime} X^{**}(X^{**\prime}X^{**})^{-1}X^{**\prime}Z_g^{**} \right\}^{-1} \\ & \times \left\{ Z_g^{**\prime}X^{**}(X^{**\prime}X^{**})^{-1} \begin{pmatrix} X_1^{**\prime}\underline{w}_g \\ X_2^{**\prime}\underline{w}_g - \hat{C}_g \end{pmatrix} \right\}, \end{aligned}$$

where $X^{**} = X\tilde{M}_x = (X_1^{**}, X_2^{**})$, $X_1^{**} = (\triangledown W_{-1}, \ldots, \triangledown W_{-p+1})$, and $X_2^{**} = W_{-p}$. Just like $(Z_{g1}^{**}, Z_{g2}^{**})$, $(X_1^{**}, X_2^{**})$ are well defined and observable.

**Theorem 5.2.** *Under assumptions A1-A4, 3.1 and 3.2, the modified 2SLS estimator $\hat{\tilde{\delta}}_{g,m2SLS}^* = D_g^{-1}\hat{\tilde{\delta}}_{g,m2SLS}^{**}$ is consistent. Furthermore*

$$(5.8) \qquad \sqrt{T}(\hat{\tilde{\delta}}_{g1,m2SLS}^* - \delta_{g1}^*) \Longrightarrow N(\underline{0}, \sigma_g^2(M_{z_{g1}x_1}^* M_{x_1x_1}^{*-1} M_{x_1z_{g1}}^*)^{-1})$$

*and is independent of*

$$(5.9) \qquad \begin{aligned} T(\hat{\tilde{\delta}}_{g2,m2SLS}^* - \delta_{g2}^*) \Longrightarrow & (M_{z_{g2}x_2}^* M_{x_2x_2}^{*-1} M_{x_2z_{g2}}^*)^{-1} \\ & \cdot M_{z_{g2}x_2}^* M_{x_2x_2}^{*-1} \int B_{x_2^*} dB_{\epsilon_g \cdot x_2^*}, \end{aligned}$$

*which is a mixed normal of the form*

$$(5.10) \qquad \int_{M_{x_2x_2}^* > 0} N(\underline{0}, \sigma_{g \cdot \triangledown x_2^*}^2 (M_{z_{g2}x_2}^* M_{x_2x_2}^{*-1} M_{x_2z_{g2}}^*)^{-1}) dP(M_{x_2x_2}^*).$$

*where $\sigma_{g \cdot \triangledown x_2^*}^2 = \sigma_g^2 - L_{\epsilon_g \triangledown x_2^*} L_{\triangledown x_2^* \triangledown x_2^*} L_{\triangledown x_2^* \epsilon_g}$.*

The modified 2SLS estimator of $\underline{\delta}_g$ can be obtained as

$$(5.11) \qquad \hat{\underline{\delta}}_{g,m2SLS} = \tilde{M}_g \hat{\underline{\delta}}_{g,m2SLS}^{**} = \tilde{M}_g D_g \hat{\underline{\delta}}_{g,m2SLS}^*,$$

where $\tilde{M}_g$ is a known matrix but in general, not $D_g$. However, although the modified 2SLS estimator of $\underline{\delta}_g^*$ is either asymptotically normal or mixed normal, the Wald type test statistic

$$(5.12) \qquad \frac{1}{\sigma_g^2}(P\hat{\underline{\delta}}_{g,m2SLS} - \underline{c})'\{P[Z_g'X(X'X)^{-1}X'Z_g]P'\}^{-1}(P\hat{\underline{\delta}}_{g,m2SLS} - \underline{c})$$

does not always have the asymptotic chi-square distribution under the null hypothesis $P\underline{\delta}_g = \underline{c}$, where $P$ is a known $k \times g_\Delta$ matrix of rank $k$. To see this, rewrite (5.12) in terms of $\hat{\underline{\delta}}_{g,m2SLS}^*$

$$(5.13) \qquad \begin{aligned} &\frac{1}{\sigma_g^2}(P^*H_g\hat{\underline{\delta}}_{g,m2SLS}^* - \underline{c})' \left\{ P^*H_g[Z_g^{*'}X^*(X^{*'}X^*)^{-1}X^{*'}Z_g^*]H_g'P^{*'} \right\} \\ &\times (P^*H_g\hat{\underline{\delta}}_{g,m2SLS}^* - \underline{c}), \end{aligned}$$

where $P^* = P\tilde{M}_g D_g H_g^{-1}$ and $H_g = \begin{bmatrix} T^{-1/2}I_{lg} & 0 \\ 0 & T^{-1}I_{b_g} \end{bmatrix}$. The null hypothesis becomes $P^*H_g\delta_g^* = \underline{c}$. Notice that the asymptotic covariance matrix of $H_g\hat{\underline{\delta}}_{g,m2SLS}^*$ converges to

$$\begin{pmatrix} \sigma_g^2(M_{z_{g1}^*x_1}^*M_{x_1x_1}^{*-1}M_{x_1z_{g1}}^*)^{-1} & \underline{0} \\ \underline{0} & \sigma_{g\cdot\nabla x_2^*}^2(M_{z_{g2}x_2}^*M_{x_2x_2}^{*-1}M_{x_2z_{g2}}^*)^{-1} \end{pmatrix},$$

while $H_g[Z_g^*X^*(X^{*'}X^*)^{-1}X^{*'}Z_g^*]H_g'$ in (5.13) converges to

$$(5.14) \qquad \sigma_g^2 \begin{pmatrix} (M_{z_{g1}x_1}^*M_{x_1x_1}^{*-1}M_{x_1z_{g1}}^*)^{-1} & \underline{0} \\ \underline{0} & (M_{z_{g2}x_2}^*M_{x_2x_2}^{*-1}M_{x_2z_{g2}}^*)^{-1} \end{pmatrix}.$$

Wald statistic (5.12) (or equivalently (5.13)) is asymptotically chi-square distributed with $k$ degrees of freedom if and only if $P\hat{\underline{\delta}}_{g,m2SLS}$ (or equivalently $P^*H_g\hat{\underline{\delta}}_{g,m2SLS}^*$) in the hypothesis does not involve the $T$-consistent component $\hat{\underline{\delta}}_{g2,m2SLS}^*$. Otherwise, $H_g[Z_g^{*'}X^*(X^{*'}X^*)^{-1}X^{*'}Z_g^{*'}]H_g'$ would overestimate the asymptotic covariance matrix of $H_g\hat{\underline{\delta}}_{g,m2SLS}^*$ because $\sigma_{g\cdot\nabla x_2^*}^2 \le \sigma_g^2$ for the submatrix corresponding to $\underline{x}_2^*$ and $\underline{z}_{g2}^*$. In general, the test statistic (5.12) is a conservative test, with its asymptotic distribution a weighted sum of $k$ independent $\chi_1^2$ variables with weights between 0 and 1.

The construction of the modified 2SLS estimator requires nonparametric estimation of the long-run covariance matrix and the one-sided long-run covariance matrix. It is well known that kernel estimator and hence the finite sample performance of the modified 2SLS estimator could be affected substantially by the choice of the bandwidth parameter. In addition, since we can not approximate the asymptotic covariance matrix of the modified 2SLS estimator properly, Wald test statistics based on the modified 2SLS estimator using the formula (5.12) may not be chi-square distributed and critical values that are based on chi-square distributions can be used for conservative tests only. However, as noted by Toda and Yamamoto [55], if we augment the order of a $p$-th order autoregressive process by the maximum order of integration then the miscentering and skewness of the limiting distribution of the least squares estimator will be concentrated on the coefficient matrices associated with the augmented lagged vectors which are known *a priori* to be zero, therefore can be ignored. Standard inference procedure can still be applied to the coefficients of the first $p$ coefficient matrices. Hsiao and Wang [28] follow this idea by proposing a lag augmented 2SLS.

The $p$-th order structural VAR (2.1) can be written as a $(p+1)$-th order structural VAR,

$$(5.15) \qquad A_0 \underline{w}_t + A_1 \underline{w}_{t-1} + \cdots + A_p \underline{w}_{t-p} + A_{p+1} \underline{w}_{t-p-1} = \underline{\epsilon}_t,$$

where $A_{p+1} \equiv \underline{0}$. Transforming (5.15) into an error-correction form, we have

$$(5.16) \qquad \sum_{j=0}^{p} A_j^* \bigtriangledown \underline{w}_{t-j} + A_{p+1}^* \underline{w}_{t-p-1} = \underline{\epsilon}_t,$$

where $A_j^* = \sum_{\ell=0}^{j} A_\ell, j = 0, 1, \ldots, p$ and $A_{p+1}^* = A_p^*$. It follows that $A = [A_0, \ldots, A_p] = [A_0^*, \ldots, A_p^*] \tilde{M}^{-1}$.

Let the $g$-th equation of (5.15) be written as

$$(5.17) \qquad \underline{w}_g = Z_g^A \underline{\delta}_g^A + \underline{\epsilon}_g,$$

where $Z_g^A = (Z_g, \tilde{\underline{w}}_{g,-(p+1)}), \underline{\delta}_g^A = (\underline{\delta}'_g, -\underline{a}'_{g,p+1})'$ with $\tilde{\underline{w}}_{g,-(p+1)}$ denoting the $T \times g_\Delta$ vector of included $\tilde{\underline{w}}_{gt}$ lagged by $(p+1)$ periods and $\underline{a}_{g,p+1}$ is the $g$-th row of $A_{p+1}$ excluding those elements subject to exclusion restrictions. Just like (4.1), there exists a nonsingular transformation matrix $M_g^A$ that transforms $Z_g^A$ into $Z_g^{*A} = Z_g^A M_g^A = (Z_{g1}^{*A}, Z_{g2}^{*A})$, and $\underline{\delta}_g^{*A} = (M_g^A)^{-1} \underline{\delta}_g^A = (\underline{\delta}_{g1}^{*A'}, \underline{\delta}_{g2}^{*A'})'$ where $Z_{g1}^{*A} = (\bigtriangledown Z_g, \tilde{W}_{g,-(p+1)} \underline{\pi}_g)$ is stationary and $Z_{g2}^{*A} = \tilde{W}_{g2,-(p+1)}$ consists of $T$ observed $b_g$ linearly independent $I(1)$ variables, $\tilde{\underline{w}}_{g2,t-(p+1)}$. Rewrite (5.17) in terms of the transformed variables,

$$(5.18) \qquad \underline{w}_g = Z_g^A M_g^A (M_g^A)^{-1} \underline{\delta}_g^A + \underline{\epsilon}_g = \begin{pmatrix} Z_{g1}^{*A} & Z_{g2}^{*A} \end{pmatrix} \begin{pmatrix} \underline{\delta}_{g1}^{*A} \\ \underline{\delta}_{g2}^{*A} \end{pmatrix} + \underline{\epsilon}_g$$

Let $X^A = (X, W_{-(p+1)})$. The 2SLS estimator of (5.17) is defined as

$$(5.19) \qquad \hat{\underline{\delta}}_{g,2SLS}^A = [Z_g^{A'} X^A (X^{A'} X^A)^{-1} X^{A'} Z_g^A]^{-1} [Z_g^{A'} X^A (X^{A'} X^A)^{-1} X^{A'} \underline{w}_g].$$

The LA2SLS of (4.1) is defined as

$$(5.20) \qquad \hat{\underline{\delta}}_{g,LA2SLS} = Q_g^A \hat{\underline{\delta}}_{g,2SLS}^A,$$

where $Q_g^A = (I_{(p+1)g_\Delta - 1}, \underline{0}_{g_\Delta})$, where $\underline{0}_{g_\Delta}$ denotes a $[(p+1)g_\Delta - 1] \times g_\Delta$ matrix of zeros. Since $\hat{\underline{\delta}}_{g,2SLS}^A = M_g^A \hat{\underline{\delta}}_{g,2SLS}^{*A}$, we have

$$(5.21) \qquad \begin{aligned} \hat{\underline{\delta}}_{g,LA2SLS} &= Q_g^A M_g^A \hat{\underline{\delta}}_{g,2SLS}^{*A} \\ &= (\tilde{M}_g, \underline{0}_{g_\Delta}) \hat{\underline{\delta}}_{g,2SLS}^{*A} \\ &= (\tilde{M}_g, \underline{0}_g) \hat{\underline{\delta}}_{g1,2SLS}^{*A}, \end{aligned}$$

where $\tilde{M}_g$ is a $[(p+1)g_\Delta - 1] \times [(p+1)g_\Delta - 1]$ matrix of the form,[8]

$$(5.22) \qquad \tilde{M}_g = \begin{pmatrix} I_{g_\Delta - 1} & \underline{0} & \dots \dots \dots \dots & \underline{0} \\ \begin{pmatrix} -I_{g_\Delta - 1} \\ \underline{0}' \end{pmatrix} & I_{g_\Delta} & \dots \dots \dots \dots & \dots \\ \dots & -I_{g_\Delta} & I_{g_\Delta} \dots \dots \dots \dots & \dots \\ \dots & \dots & \dots \dots \dots \dots & I_{g_\Delta} \quad \underline{0} \\ \dots & \dots & \dots \dots \dots \dots & -I_{g_\Delta} \quad I_{g_\Delta} \end{pmatrix},$$

with $w_{gt}$ being put as the last element of $\tilde{\underline{w}}_{gt}$, $I_{g_\Delta}$ denoting the identity matrix of the dimension of included variables in the $g$-th equation, and $\underline{0}_g$ is a $[(p+1)g_\Delta - 1] \times r_g$ matrix with $r_g$ denoting the number of cointegrating relations among $\tilde{\underline{w}}_{gt}$ such that $\underline{w}'_{gt}\tilde{\underline{\pi}}_g$ is $I(0)$. Then $\underline{\delta}_g = (\tilde{M}_g, \underline{0}_g)\underline{\delta}^{*A}_{g1}$.

Since

$$(5.23) \qquad \sqrt{T}(\hat{\underline{\delta}}^{*A}_{g1,2SLS} - \underline{\delta}^{*A}_{g1}) \longrightarrow N[\underline{0}, \sigma^2_g(M^{A*}_{z_{g1}x_1} M^{A*-1}_{x_1x_1} M^{A*}_{x_1z_{g1}})^{-1}],$$

where $M^{A*}_{z_{g1}x_1} = \text{plim } \frac{1}{T} Z^{*A'}_{g1} X^{*A}_1, M^{A*}_{x_1x_1} = \text{plim } \frac{1}{T} X^{*A'}_1 X^{*A}_1$, with $X^{*A}_1 = (\bigtriangledown X, W_{-(p+1)}\underline{d})$ being the $T \times (mp + r)$ linearly independent $I(0)$ variables. It follows that

**Theorem 5.3.** *The LA2SLS of $\delta_g$ is consistent and*

$$(5.14) \qquad \sqrt{T}(\hat{\delta}_{g,LA2SLS} - \delta_g)$$
$$\Longrightarrow N \left\{ \underline{0}, \sigma^2_g(\tilde{M}_g \quad \underline{0}_g)[M^{A*}_{z_{g1}x_1} M^{A*-1}_{x_1x_1} M^{A*}_{x_1z_{g1}}]^{-1} \begin{pmatrix} \tilde{M}'_g \\ \underline{0}'_g \end{pmatrix} \right\}.$$

The LA2SLS estimators of the coefficients of the original structural VAR model (2.1) converge to the true value at the speed of $T^{1/2}$ and are asymptotically normally distributed with nonsingular covariance matrix. Therefore, Wald type test statistics based on LA2SLS estimates are asymptotically chi-square distributed. Compared to the conventional 2SLS or modified 2SLS, the LA2SLS estimator loses the T-convergence component and ignores the prior restrictions that the coefficients on $\tilde{\underline{w}}_{g,t-(p+1)}$ are zero, hence may lose some efficiency. However, since distribution of $\hat{\underline{\delta}}_g$ is a linear combination of $\hat{\underline{\delta}}^*_{g1}$ and $\hat{\underline{\delta}}^*_{g2}$ and the limiting distribution of $\hat{\underline{\delta}}_{g,LA2SLS}$ is given by the components of the slower rate of convergence, the loss of efficiency in estimating $\hat{\underline{\delta}}_g$ by LA2SLS may not be that significant, as reported in a Monte Carlo Study by Hsiao and Wang [28].

## 6. Conclusions

As demonstrated by Nelson and Plosser [41] that many economic time series are nonstationary. The advancement of nonstationary time series analysis provides a rich reportoire of analytic tools for economists to analyze how do variables respond dynamically to shocks through the decomposition a dynamic system into long-run and short-run relations and allow economists to extract common stochastic trends present in the system that provide information on the important sources of economic fluctuation (e.g. Banerjee, Dolado, Galbraith and Hendry [5], King, Plosser, Stock and Watson [33]). However nonstationarity does not invalid the main concerns of

---

[8]For ease of notation, we assume all the included variables appear with the same lag order.

Cowles Commission structural approach — identification and simultaneity bias. As shown by Hsiao [26], whether the data is stationary or nonstationary, the same rank condition holds for the identification of an equation in a system. Ignoring the correlations between the regressors and the errors of the equation that arise from the joint dependency of economic variables can lead to severe bias in the least squares estimator even though the regressors are I(1) (Hsiao [21], also see the Monte Carlo study by Hsiao and Wang [28]). Instrumental variable methods have to be applied to obtain consistency.

However, nonstationarity does raise the issue of statistical inference. Standard instrumental variable method can lead to estimators that have non-normal asymptotic distributions and are asymptotically biased and skewed. If there exists prior knowledge to dichotomize the set of variables into joint dependent and exogenous variables and the nonstationarity in the dependent variables is driven by the nonstationarity in the exogenous variables through cointegration relations, standard 2SLS developed for the stationary data can also be used for the analysis of nonstationary data (Hsiao [21, 22]). Wald test statistics for the null are asymptotically chi-square distributed. There is no inference issue. On the other hand, if all the variables are treated as joint dependent as in the time series context, although 2SLS is consistent, the limiting distribution is subject to miscentering and skewness associated with the unit root distribution. Modified or lag order augmented 2SLS will have to be used to ensure valid inference. The modified 2SLS is asymptotically more efficient. However, it also suffers more size distortion in finite sample. On the other hand, the lag order augmented 2SLS does not suffer much efficiency loss, at least in a small scale SVAR model (e.g. Hsiao and Wang [28]), and chi-square distribution is a good approximation for the test statistic.

All above discussions were based on the assumption that no knowledge of cointegration or direction of nonstationarity is known *a priori*. If such information is available, (e.g. King, Plosser, Stock and Watson [33]) estimators incorporating the knowledge of the rank of cointegration presumably will not only lead to efficient estimators of structural form parameters, but also avoid the inference issues arising from the matrix unit roots distrubutions in the system. Unfortunately, structural form estimation methods incorporating reduced rank restrictions appear to be fairly complicated.

The focus of this review is to take a SVAR model as a maintained hypothesis, search for better estimators and understand their properties. We have not looked at the issues of modeling strategy. There is a vast literature on the interactions between structural and non-structural time series analysis to uncover the data-generation process, including testing, estimation, model-combining and prediction (e.g. Hendry and Ericsson [16], Hendry and Krolzig [17], King, Plosser, Stock and Watson [33], Zellner and Palm [61]).

## Acknowledgments

I would like to thank Peter Robinson, Arnold Zellner and a referee for helpful comments.